\documentstyle[11pt]{article}
\def\b#1{{\bf #1}}
\def\i#1{{\it #1}}

\begin{document}
\title{Identities for Fibonacci and Lucas Polynomials derived
from a book of Gould}
\author{Mario Catalani\\
Department of Economics, University of Torino\\ Via Po 53, 10124 Torino, Italy\\
mario.catalani@unito.it}
\date{}
\maketitle
\begin{abstract}
\small{This note is dedicated to Professor Gould. The aim is to show
how the identities in his book "Combinatorial Identities" can be used
to obtain identities
for Fibonacci and
Lucas polynomials. In turn these identities allow to derive a wealth
of numerical identities for Fibonacci and Lucas numbers.}
\end{abstract}

\section{Introduction}
The book of Gould \cite{gould} is an almost endless source of applications
in many fields: in this note we give examples showing how identities for
Fibonacci and Lucas polynomials can be derived. As a consequence many
numerical identities can be obtained choosing appropriate numerical values
for the variables $x$ and $y$.

\bigskip
\noindent
We define bivariate Fibonacci polynomials as
$$F_n(x,\,y)=xF_{n-1}(x,\,y)+yF_{n-2}(x,\,y),\qquad F_0(x,\,y)=0,\,
F_1(x,\,y)=1,$$ and
bivariate Lucas polynomials as
$$L_n(x,\,y)=xL_{n-1}(x,\,y)+yL_{n-2}(x,\,y),\qquad L_0(x,\,y)=2,\,
L_1(x,\,y)=x.$$
We assume $x\not= 0,\,y\not= 0,\,x^2+4y\not= 0$.

\noindent
The roots of the characteristic equation are
$$\alpha\equiv\alpha(x,\,y)={x+\sqrt{x^2+4y}\over 2},\quad\quad
\beta\equiv\beta(x,\,y)={x-\sqrt{x^2+4y}\over 2}.$$
We have $\alpha+\beta=x,\,\alpha\beta=-y$ and $\alpha-\beta=\sqrt{x^2+4y}$.
The Binet's forms are
$$F_n(x,\,y)={\alpha^n-\beta^n\over \alpha-\beta},
\quad\quad L_n(x,\,y)=\alpha^n+\beta^n.$$
The generating function of the Fibonacci polynomials is
$${t\over 1-xt-yt^2},$$
that of the Lucas polynomials is
$${2-xt\over 1-xt-yt^2}.$$
Many basic facts concerning these kinds of polynomials can be found in
\cite{mario7}, \cite{mario8}.

\section{Examples}
All the references are to the book \cite{gould}.
\bigskip

\noindent
\b{1.} Identity 1.64 says
$$\sum_{k=0}^{\left\lfloor{n\over 2}\right\rfloor}{n-k\choose k}
{1\over n-k}\left ({z\over 4}\right )^k=
{1\over n2^{n-1}}{u^n+v^n\over u+v},$$
with $u=1+\sqrt{z+1},\,v=1-\sqrt{z+1}$.

\noindent
Write $z=4{y\over x^2}$. Then
$$u=1+\sqrt{1+4{y\over x^2}}={x+\sqrt{x^2+4y}\over x},$$
from which $u={2\alpha\over x}$. In the same way
$v={x-\sqrt{x^2+4y}\over x}= {2\beta\over x}.$
The RHS of the previous identity becomes
$${L_n(x,\,y)\over nx^n},$$
so that we get the identity
$$
\sum_{k=0}^{\left\lfloor{n\over 2}\right\rfloor}{n-k\choose k}
{n\over n-k}x^{n-2k}y^k=L_n(x,\,y).$$

\bigskip
\noindent
\b{2.} Identity 1.38 says
$$
\sum_{k=0}^{\left\lfloor{n\over 2}\right\rfloor}{n\choose 2k}
{z^{2k}\over 2k+1}={(1+z)^{n+1}-(1-z)^{n+1}\over 2(n+1)z}.$$
Write $z={\alpha-\beta\over x}$. Then $1+z={2\alpha\over x},\,
1-z={2\beta\over x}$. The RHS of the identity becomes
$${2^n\over (n+1)x^n}F_{n+1}(x,\,y),$$
so that we get the identity
$$
\sum_{k=0}^{\left\lfloor{n\over 2}\right\rfloor}{n\choose 2k}
{(x^2+4y)^kx^{n-2k}\over 2k+1}={2^n\over (n+1)}F_{n+1}(x,\,y).$$

\bigskip
\noindent
\b{3.} Identity 1.39 says
$$
\sum_{k=0}^{\left\lfloor{n-1\over 2}\right\rfloor}{n\choose 2k+1}
{z^{2k}\over k+1}={(1+z)^{n+1}+(1-z)^{n+1}-2\over (n+1)z^2}.$$
Multiply both sides by $z^2$ and again
write $z={\alpha-\beta\over x}$.
The identity becomes
$$
\sum_{k=0}^{\left\lfloor{n-1\over 2}\right\rfloor}{n\choose 2k+1}
{(x^2+4y)^{k+1}\over (k+1)x^{2k+2}}={2^{n+1}L_{n+1}(x,\,y)-2x^{n+1}\over
(n+1)x^{n+1}},$$
that is
$$
(n+1)\sum_{k=0}^{\left\lfloor{n-1\over 2}\right\rfloor}{n\choose 2k+1}
{(x^2+4y)^{k+1}x^{n-2k-1}\over k+1}=2^{n+1}L_{n+1}(x,\,y)-2x^{n+1}.$$

\bigskip
\noindent
\b{4.} Identity 1.61 says
$$
\sum_{k=0}^{\left\lfloor{n\over 2}\right\rfloor}{n-k\choose k}
2^{n-2k}y^k={\left (1+\sqrt{1+y}\right )^{n+1}-
\left (1-\sqrt{1+y}\right )^{n+1}\over 2\sqrt{1+y}}.$$
Consider the polynomials $F_n(2,\,y)$. The roots of the characteristic
equation are
$$\alpha= 1+\sqrt{1+y},\qquad \beta=1-\sqrt{1+y}.$$
It follows
$$
\sum_{k=0}^{\left\lfloor{n\over 2}\right\rfloor}{n-k\choose k}
2^{n-2k}y^k={\alpha^{n+1}-\beta^{n+1}\over \alpha-\beta}=F_{n+1}(2,\,y).$$
Recall that $F_n(2,\,1)$ are the Pell numbers $P_n$.

\bigskip
\noindent
\b{5a.} Identity 1.87 says
$$
\sum_{k=0}^{\left\lfloor{n\over 2}\right\rfloor}{n\choose 2k}z^k
={(1+\sqrt{z})^n+(1-\sqrt{z})^n\over 2}.$$
With $z=1+4{y\over x^2}$ we have
$$1+\sqrt{z}={x+\sqrt{x^2+4y}\over x}={2\alpha\over x},$$
$$1-\sqrt{z}={x-\sqrt{x^2+4y}\over x}={2\beta\over x}.$$
Then
$$\sum_{k=0}^{\left\lfloor{n\over 2}\right\rfloor}{n\choose 2k}
\left ({x^2+4y\over x^2}\right )^k={2^{n-1}\over x^n}L_n(x,\,y),$$
that is
$$\sum_{k=0}^{\left\lfloor{n\over 2}\right\rfloor}{n\choose 2k}
(x^2+4y)^kx^{n-2k}=2^{n-1}L_n(x,\,y).$$

\bigskip
\noindent
\b{5b.} In the same Identity 1.87 replace $z$ by ${\alpha^2\over x^2}$
so that
$$1+\sqrt{z}={x+\alpha\over x},\qquad 1-\sqrt{z}={\beta\over x},$$
and then by
${\beta^2\over x^2}$
so that now
$$1+\sqrt{z}={x+\beta\over x},\qquad 1-\sqrt{z}={\alpha\over x}.$$
After summation, the identity becomes
$$\sum_{k=0}^{\left\lfloor{n\over 2}\right\rfloor}{n\choose 2k}
{L_{2k}(x,\,y)\over x^{2k}}={1\over 2x^n}L_n(x,\,y)+
{1\over 2}\left [\left ({x+\alpha\over x}\right )^n+
\left ({x+\beta\over x}\right )^n\right ].$$
Let
$${\tilde \alpha}=x+\alpha={3x+\sqrt{x^2+4y}\over 2},$$
$${\tilde \beta}=x+\beta={3x-\sqrt{x^2+4y}\over 2}.$$
Then ${\tilde \alpha}$ and ${\tilde \beta}$ are the roots of the
characteristic equation of $L_n(3x,\,y-2x^2)$ as it can be checked
easily. Then
$$\sum_{k=0}^{\left\lfloor{n\over 2}\right\rfloor}{n\choose 2k}
{L_{2k}(x,\,y)\over x^{2k}}={1\over 2x^n}\left [L_n(x,\,y)+
L_n(3x,\,y-2x^2)\right ],$$
that is
$$2\sum_{k=0}^{\left\lfloor{n\over 2}\right\rfloor}{n\choose 2k}
L_{2k}(x,\,y)x^{n-2k}=L_n(x,\,y)+
L_n(3x,\,y-2x^2).$$

\bigskip
\noindent
\b{6a.} Identity 1.95 says
$$
\sum_{k=0}^{\left\lfloor{n-1\over 2}\right\rfloor}{n\choose 2k+1}z^k
={(1+\sqrt{z})^n-(1-\sqrt{z})^n\over 2\sqrt{z}}.$$
As in Example 5a, write
$z=1+4{y\over x^2}$. Then
$$\sum_{k=0}^{\left\lfloor{n-1\over 2}\right\rfloor}{n\choose 2k+1}
(x^2+4y)^kx^{n-2k-1}=2^{n-1}F_n(x,\,y).$$

\bigskip
\noindent
\b{6b.} In the same Identity 1.95 replace, as in Example 5b,
$z$ a first time by ${\alpha^2\over x^2}$
and then by
${\beta^2\over x^2}$.
Subtract the resulting identities and divide by $\alpha-\beta$.
The LHS is now
$$\sum_{k=0}^{\left\lfloor{n-1\over 2}\right\rfloor}{n\choose 2k+1}
{F_{2k}(x,\,y)\over x^{2k}}.$$
The RHS after some simple manipulations becomes
$${1\over 2x^{n-1}(\alpha-\beta)}\left [{(x+\alpha)^n-\beta^n\over
\alpha}-
{(x+\beta)^n-\alpha^n\over
\beta}\right ],$$
that is
$${1\over 2x^{n-1}(-y)(\alpha-\beta)}\left [(\alpha-\beta)F_{n+1}(x,\,y)
+\beta (x+\alpha)^n-\alpha (x+\beta)^n\right ].$$
As in Example 5b, let
$${\tilde \alpha}=x+\alpha={3x+\sqrt{x^2+4y}\over 2},$$
$${\tilde \beta}=x+\beta={3x-\sqrt{x^2+4y}\over 2}.$$
Again ${\tilde \alpha}$ and ${\tilde \beta}$ are the roots of the
characteristic equation of $L_n(3x,\,y-2x^2)$
or $F_n(3x,\,y-2x^2)$. Then ${\tilde \alpha}-{\tilde \beta}=\alpha
-\beta$ and ${\tilde \alpha}{\tilde \beta}=-y+2x^2$. It follows
\begin{eqnarray*}
\beta (x+\alpha)^n&=& ({\tilde \beta}-x){\tilde \alpha}^n\\
&=& {\tilde \beta}{\tilde \alpha}^n-x{\tilde \alpha}^n\\
&=& {\tilde \beta}{\tilde \alpha}{\tilde \alpha}^{n-1}-x{\tilde \alpha}^n\\
&=& (-y+2x^2){\tilde \alpha}^{n-1}-x{\tilde \alpha}^n.
\end{eqnarray*}
In the same way
$$\alpha (x+\beta)^n=(-y+2x^2){\tilde \beta}^{n-1}-x{\tilde \beta}^n.$$
So the RHS now becomes
$${1\over 2x^{n-1}(-y)}\left [F_{n+1}(x,\,y)-(y-2x^2)
F_{n-1}(3x,\,y-2x^2)-xF_n(3x,\,y-2x^2)\right ].$$
Finally we get the identity
\begin{eqnarray*}
&&2y\sum_{k=0}^{\left\lfloor{n-1\over 2}\right\rfloor}{n\choose 2k+1}
F_{2k}(x,\,y)x^{n-2k-1}\\ &&\quad\quad =
-\left [
F_{n+1}(x,\,y)-(y-2x^2)
F_{n-1}(3x,\,y-2x^2)-xF_n(3x,\,y-2x^2)\right ].
\end{eqnarray*}

\bigskip
\noindent
\b{Remark.} $F_n(1,\,1)$ and $L_n(1,\,1)$ are, respectively, the Fibonacci
and Lucas numbers. In the above Examples, with the exception of 4, 5b, 6b,
all the identities can be reduced to identities involving Fibonacci
and Lucas numbers by setting $x=L_k(1,\,1),\,y=(-1)^{k+1}$. This is due
to a consequence of a result in \cite{mario8}:
$$F_n\left (L_k(x,\,y),\,(-1)^{k+1}y^k\right )=
{F_{kn}(x,\,y)\over F_k(x,\,y)},$$
$$L_n\left (L_k(x,\,y),\,(-1)^{k+1}y^k\right )=
L_{kn}(x,\,y).$$

\end{document}